\documentclass[12pt]{article}
\usepackage{amssymb}
\usepackage{epsfig}

\newtheorem{prop}{Proposition}[section]
\newtheorem{thm}[prop]{Theorem}
\newtheorem{lemma}[prop]{Lemma}
\newtheorem{cor}[prop]{Corollary}

\newenvironment{proof}{\noindent{\em Proof.}~}{{\hfill $\square$}\par\medskip}

\def\Z{\mathbb{Z}}
\def\R{\mathbb{R}}
\def\C{\mathbb{C}}
\def\Q{\mathbb{Q}}
\def\G{\Gamma}
\def\<{\langle}
\def\>{\rangle}
\def\sm{\smallsetminus}
\def\e{\varepsilon}

\begin{document}
\title{Winding numbers and $SU(2)$-representations of knot groups}

\author{Dylan Bowden and James Howie\\
Department of Mathematics and the\\
Maxwell Institute for Mathematical Sciences\\
Heriot-Watt University\\
Edinburgh EH14 4AS\\
UK}

\maketitle

\begin{abstract}
Given an abelian group $A$ and a Lie group $G$, we construct a
bilinear pairing from $A\times\pi_1({\mathcal R})$ to $\pi_1(G)$,
where $\mathcal R$ is a subvariety of
the variety of representations $A\to G$.

In the case where $A$ is the peripheral subgroup of a torus or two-bridge
knot group, $G=S^1$ and $\mathcal R$ is a certain variety of representations
arising from suitable $SU(2)$-representations of the knot group, we show
that this pairing is not identically zero. We discuss the consequences
of this result for the $SU(2)$-representations of fundamental
groups of manifolds obtained by Dehn surgery on such knots.
\end{abstract}

\section{Introduction}

The real algebraic variety of representations from a $3$-manifold
group $\pi_1(M)$ to $SU(2)$ or $SO(3)$ has long been a subject of
interest, giving rise as it does to useful invariants such as the
Casson invariant and the $A$-polynomial \cite{CCGLS}.

In the case where $\partial M$ is a torus -- in particular, where $M$
is the exterior of a knot in $S^3$ -- there is a particular interest in finding representations which vanish on a given slope $\alpha\in\Q\cup\{\infty\}$ on $\partial M$,
and hence give rise to a representation of $\pi_1(M(\alpha))$, where
$M(\alpha)$ is the manifold obtained from $M$ by Dehn filling along $\alpha$.

A description of the character variety in the case of a $2$-bridge
knot is given by Burde in \cite{B}.  For twist knots, a more detailed
description is given by Uygur and Azcan in \cite{UA}.

Burde \cite{B} used this description to show that nontrivial representations $\pi_1(M(+1))\to SU(2)$
exist for any nontrivial $2$-bridge knot exterior $M$, and
deduced the Property P Conjecture for $2$-bridge knots.  More recently,
Kronheimer and Mrowka \cite{KM1} proved the Property P Conjecture in full by
showing that nontrivial representations $\pi_1(M(+1))\to SO(3)$ exist for 
an arbitrary nontrivial knot exterior $M$.

In another article \cite{KM2}, the same authors proved that there
is an irreducible representation $\pi_1(M(r))\to SU(2)$ 
(that is, a representation with nonabelian image), for
any nontrivial knot exterior $M$ and any slope $r\in\Q$ such that
$|r|\le 2$.  One consequence of this (see \cite{BZ,DG})
is that every nontrivial knot has a nontrivial $A$-polynomial.

\medskip
In the present note, we construct a bilinear pairing 
$\pi_1({\mathcal C})\times\pi_1(\partial M)\to\Z$ for suitable subsets
$\mathcal C$ of the variety $\mathcal R$ of representations 
$\pi_1M\to SU(2)$, and apply it to Burde's description \cite{B} of
$\mathcal R$ in the case of $2$-bridge knots, to show that the 
restriction $|r|\le 2$ in \cite{KM2} can be weakened in this case:

\begin{thm}\label{twobridge}
Let $M$ be the exterior of a nontrivial $2$-bridge knot in $S^3$
which is not a torus knot,
and let $\alpha$ be any non-meridian slope in $\partial M$.  Then
there exists an irreducible representation $\pi_1(M(\alpha))\to SU(2)$.
\end{thm}

Since there are many examples of lens spaces obtainable by Dehn surgery
on nontrivial knots, it is clear that the above theorem cannot possibly
extend from $2$-bridge knots to arbitrary knots.  However, by varying
the subset $\mathcal C$ of the representation in our construction,
we can adapt the technique to consider also reducible representations.
 
As an example, we prove the following result for torus knots.

\begin{thm}\label{torus}
Let $X$ be the exterior of the $(p,q)$ torus knot, where $1<p<q$,
and $X(\alpha)$
the manifold obtained from $X$ by Dehn filling along a non-meridian
slope $\alpha\in\Q$.  Then
\begin{enumerate}
\item[\rm{(1)}] if $\alpha=pq$ and $p>2$, then $\pi_1(X(\alpha))$ admits
an irreducible representation to $SU(2)$;
\item[\rm{(2)}] if $\alpha=pq$ and $p=2$, then $\pi_1(X(\alpha))$ admits
no irreducible representation to $SU(2)$, but admits
a representation to $SO(3)$ with nonabelian image;
\item[\rm{(3)}] if $\alpha=pq\pm\frac1n$ for some positive integer
$n$, then every representation from $\pi_1(X(\alpha))$ to $SO(3)$
has abelian image;
\item[\rm{(4)}] for any other value of $\alpha$, $\pi_1(X(\alpha))$ admits
an irreducible representation to $SU(2)$.
\end{enumerate}
\end{thm}

Results of \cite{M} indicate that this result is in a sense
best possible: for example, in
Case (3) the Dehn surgery manifold $X(\alpha)$ is a lens space.

The paper is organised as follows.  In Section \ref{charvar}
below we recall some basic properties of the $SU(2)$
representation and character varieties of a knot group.
In Section \ref{wind}
we describe our bilinear pairing, in a fairly general context.
We then apply this in Sections \ref{2bk} and \ref{tk}  to prove Theorems
\ref{twobridge} and \ref{torus}  respectively.

\medskip\noindent{\bf Acknowledgement} We are grateful to Ben Klaff
for helpful conversations about this work.

\section{The SU(2) representation and character varieties}\label{charvar}

If $\Gamma$ is any finitely presented group,
and $G$ is a (real) algebraic matrix group, then the set of
representations $\Gamma\to G$ forms a real affine algebraic variety
$\mathcal{R}$ on which $G$ acts by conjugation, giving rise to a
quotient {\em character variety} $\mathcal{X}$.

For the purposes of the present paper, $\Gamma$ will always be a knot group,
and $G=SU(2)$.  In this case $\mathcal{R}$ is naturally expressed as a union of
two closed $SU(2)$-invariant subvarieties $\mathcal{R}_{red}\cup\mathcal{R}_{irr}$,
and hence also $\mathcal{X}$ is a union of subvarieties $\mathcal{X}_{red}\cup\mathcal{X}_{irr}$.  Here $\mathcal{R}_{red}$ denotes the variety of
{\em reducible} representations $\rho:\Gamma\to SU(2)$, in other words
those for which the resulting $\Gamma$-module $\C^2$ splits as a direct sum of
two $1$-dimensional modules.  This happens precisely when the image of 
$\rho$ is abelian, in other words when $\rho$ is induced from a representation
of $\Gamma/[\Gamma,\Gamma]\cong\Z$.  Hence $\mathcal{R}_{red}$ is canonically
homeomorphic to $SU(2)\cong S^3$.  The corresponding character subvariety
$\mathcal{X}_{red}$ is canonically homeomorphic to the closed interval $[-2,2]\subset\R$, parametrised by the trace $Tr(\rho(\mu))$, where $\rho$
is a representative of a conjugacy class of reducible representations, and
$\mu\in\Gamma$ is a fixed meridian element.
The complement of $\mathcal{R}_{red}$ in $\mathcal{R}$ is not closed, but its 
closure is a subvariety $\mathcal{R}_{irr}$ which is $SU(2)$-invariant and hence
gives rise to a closed subvariety $\mathcal{X}_{irr}$ of $\mathcal{X}$.

Now fix once and for all a meridian $\mu\in\Gamma$, and consider the following 
subset $\mathcal{C}$ of $\mathcal{R}$.  A representation $\rho:\Gamma\to SU(2)$
belongs to $\mathcal{C}$ if and only if
$$\rho(\mu)=\left(\begin{array}{cc} x+iy & 0 \\ 0 & x-iy\end{array}\right)$$
with $x,y\in\R$, $x^2+y^2=1$ and $y\ge 0$.
Note that every representation in $\mathcal{R}$ is conjugate to one in
$\mathcal{C}$, so the quotient map $\mathcal{R}\to\mathcal{X}$ restricts
to a surjection on $
\mathcal{C}$ (and to a homeomorphism $
\mathcal{C}\cap\mathcal{R}_{red}\to\mathcal{X}_{red}$).

\section{Winding numbers}\label{wind}

Let $A$ be an abelian group, $G$ a (connected)
Lie group, and ${\mathcal D}$ a subset of the
variety of representations $A\to G$.  Given any path 
$P=\{\rho_t,0\le t\le 1\}$ in ${\mathcal D}$, and any $a\in A$, we obtain a path
$P(a)=\{\rho_t(a),0\le t\le 1\}$ in $G$.

  Clearly, if $P'$ is homotopic
(rel end points) to $P$, then $P'(a)$ is homotopic (rel end points)
to $P(a)$, for any $a\in A$.  Hence we obtain a pairing
$$\nu:A\times\pi_1({\mathcal D})\to\pi_1(G),~~~~\nu(a,[P]) := [P(a)].$$

\noindent{\bf Remark} Recall that, if $f,g:[0,1]\to G$ are closed paths
in the topological group $G$, based at the identity element $1\in G$,
then $[f][g]=[f.g]=[g][f]$ in $\pi_1(G,1)$, where $f.g$ denotes the
pointwise product $t\mapsto f(t)g(t)\in G$.  This can easily be seen, for
example, from the diagram below, representing the map $[0,1]^2\to G$,
$(s,t)\mapsto f(s)g(t)$.

\begin{center}
\epsfxsize=5cm \epsfysize=4.5cm \epsfbox{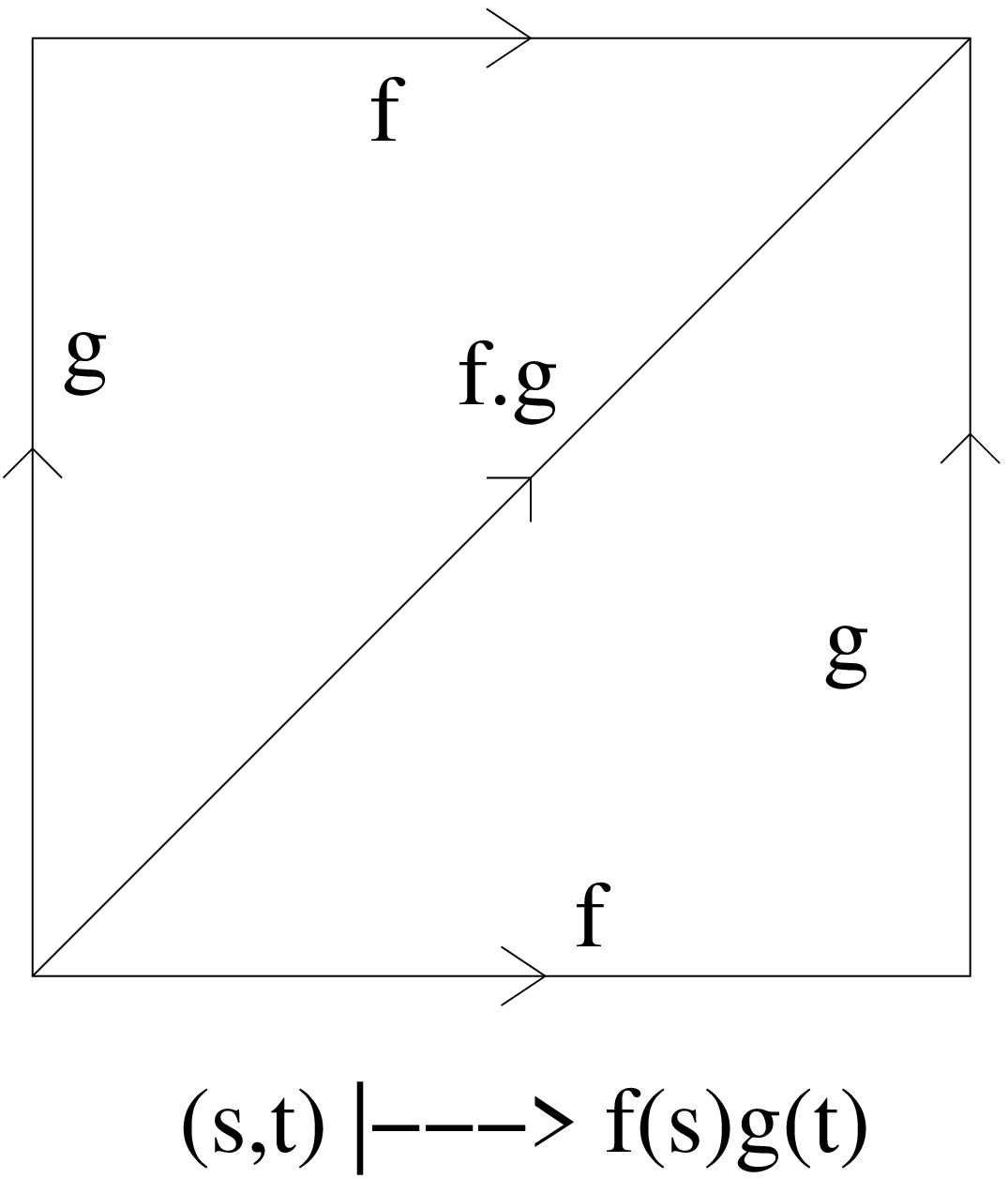}
\end{center}

In particular, $\pi_1(G,1)$ is abelian, so the above definition of
$\nu$ is unaffected by base-point choices.

\begin{prop}\label{bilinear} The pairing $\nu:(a,[P])\mapsto [P(a)]$ defined
above is bilinear.
\end{prop}

\begin{proof}
For a fixed element $a\in A$, if $P.Q$ is the concatenation of
paths $P,Q$ in ${\mathcal D}$, then $(P.Q)(a)$ is the concatenation of
$P(a),Q(a)$ in $G$, so $[P]\mapsto [P(a)]$ is a homomorphism
$\pi_1({\mathcal D})\to\pi_1(G)$.

\smallskip
For $a,b\in A$ and a fixed path $P$ in ${\mathcal D}$, we have $P_t(ab)=P_t(a)P_t(b)$
for each $t\in [0,1]$, since $P_t$ is a representation $A\to G$. By the above remark, $[P(ab)]=[P(a)][P(b)]$ in $\pi_1(G)$.
In other words, $a\mapsto [P(a)]$ is a homomorphism $A\to\pi_1(G)$.
\end{proof}

We apply Proposition \ref{bilinear} in the following restricted context.  Let $X$ be the
exterior of a nontrivial knot in $S^3$, and let 
$A=\pi_1(\partial X)\cong\Z^2$.  Let $\mu,\lambda\in A$ denote a fixed meridian and longitude respectively.

Let $G$ be the Lie group $S^1=\{z\in\C;|z|=1\}$.

The subset $\mathcal{D}$ of the variety of representations $A\to S^1$
arises as follows.  We regard $S^1$ as the subgroup of $SU(2)$ consisting
of diagonal matrices.  Recall that $\mathcal{R}$ is the variety of 
representations $\pi_1(X)\to SU(2)$, and that $\mathcal{C}$ is the
subvariety of $\mathcal{R}$ consisting of representations
$\rho:\pi_1(X)\to SU(2)$ such that $\rho(\mu)$ is diagonal, and
the imaginary part of the $(1,1)$ entry of $\rho(\mu)$ is non-negative.
Since $A$ is abelian and $\pi_1(X)$ is generated by conjugates of $\mu$,
it follows that $\rho(A)$ contains only diagonal matrices whenever
$\rho\in\mathcal{C}$.  We define $\mathcal{D}$ to be the set
of representations $A\to S^1$ that arise as restrictions of
representations in $\mathcal{C}$.

Note that $\pi_1(S)\cong\Z$, so the bilinear pairing $\nu:A\times\pi_1(\mathcal{D})\to\pi_1(S^1)$ is integer-valued.

\begin{prop}
For each $\gamma\in\pi_1(\mathcal{D})$ let $K_\gamma\subset A$ denote
the kernel of the homomorphism $A\to\Z$, $a\mapsto\nu(a,\gamma)$.
Then either $K_\gamma=A$ or $K_\gamma=\Z\mu$, the subgroup of $A$
generated by $\mu$.
\end{prop}

\begin{proof}
Certainly $\mu$ belongs to $K_\gamma$ for all
$\gamma\in\pi_1(\mathcal{D})$, since for $\rho\in\mathcal{C}$ we have $\rho(\mu)$ contained in an
open interval in $S^1$ (so the winding number of $\rho(\mu)$
as $\rho$ travels around $C$ is zero).

On the other hand, let $c=\nu(\lambda,\gamma)$.  Then by bilinearity,
for any $m,n\in\Z$ we have $\nu(m\mu+n\lambda,\gamma)=cn$.  If $cn=0$
for some $n$  then either $c=0$ or $n=0$.  In the first
case $\nu(m\mu+n\lambda,\gamma)=0$ for all $m,n$.  
In the second case, $m\mu+n\lambda=m\mu\in\Z\mu$.
\end{proof}

\begin{cor}\label{dehnfill}
If the pairing $\nu:A\times\pi_1(\mathcal{D})\to\Z$ is not uniformly vanishing, and $\alpha$ is
any non-meridian slope on $\partial X$, then $\pi_1(X(\alpha))$ admits
a nontrivial representation to $SU(2)$, where $X(\alpha)$ is the
$3$-manifold obtained from $X$ by Dehn-filling along $\alpha$.
\end{cor}

\begin{proof}
By hypothesis, $K_\gamma\ne A$ for some $\gamma\in\pi_1(\mathcal{D})$,
so $K_\gamma=\Z\mu$ by the Proposition.
Since $\alpha\notin\Z\mu$, it follows that $\nu(\alpha,\gamma)\ne 0$.
Hence the map $S^1\to S^1$ defined by
$t\mapsto\gamma_t(\alpha)$, has nonzero winding number, and hence in
particular is surjective.  Thus we may choose $t\in S^1$ such that $\gamma_t(\alpha)=1\in S^1$.  Now $\gamma_t$ is the restriction of a nontrivial representation $\sigma:\pi_1(X)\to SU(2)$, so $\sigma(\alpha)=1\in SU(2)$ and hence $\sigma$ induces a nontrivial representation
$$\tau:\pi_1(X(\alpha))=\pi_1(X)/\<\<\alpha\>\>\to SU(2).$$
\end{proof}

In practice, to find suitable closed paths in $\mathcal{D}$
we may find a closed path in $\mathcal{C}$ and project it to
$\mathcal{D}$ using the restriction map $\rho\mapsto\rho|_A$.
The next result shows that it is equally valid to work in the
character variety $\mathcal{X}$ rather than $\mathcal{C}$.

\begin{lemma}\label{char}
The restriction map $\mathcal{C}\to\mathcal{D}$, $\rho\mapsto\rho|_A$,
factors through $\mathcal{X}$.
\end{lemma}

\begin{proof}
Given $\rho,\rho'\in\mathcal{C}$ with the same image in $\mathcal{X}$,
we know that $\rho,\rho'$ are conjugate by some matrix $M\in SU(2)$.
If $\rho(\mu)\in Z(SU(2))=\{\pm I\}$, then the image of $\rho$
is central and so $\rho'=\rho$.  Otherwise, $\rho(\mu)=\rho'(\mu)$
is a diagonal matrix with non-real diagonal entries, so the conjugating
matrix $M$ must also be diagonal.  But in this case $\rho(A)$ consists
only of diagonal matrices, which therefore commute with $M$, so
the restrictions of $\rho$ and $\rho'$ to $A$ coincide.
\end{proof}

An immediate consequence of Lemma \ref{char} is that any path in
$\mathcal{R}$ between two conjugate representations gives rise to
a closed path in $\mathcal{D}$ by first projecting to $\mathcal{X}$
and then applying the restriction map $\mathcal{X}\to\mathcal{D}$.

\section{Two-bridge knots}\label{2bk}

In this section we prove the following result.

\begin{thm}\label{twobridgethm}
Let $\mathcal{R}_{irr}$ be the variety of irreducible
$SU(2)$-representations of a
two-bridge knot group $G$, and let $A$ be a peripheral subgroup
of $G$.  Then there is a closed curve $\gamma$ in $\mathcal{R}_{irr}$
such that the pairing $\nu:\pi_1(\gamma)\times A\to\Z$ is not identically zero.
\end{thm}

\begin{proof}
A two-bridge knot group $G$ has a presentation of the form

$$G=\< x, y | Wx=yW \>,$$
where $W=W(x,y)$ is a word of the form 
$x^{\e(1)}y^{\e(2)}\cdots y^{\e(2n)}$ with $\e(i)=\pm 1$ for
each $i$.  Here $x$ and $y$ are meridians.
The symmetry of the presentation ensures that $xW^*=W^*y$
in $G$, where $W^*(x,y):=W(y,x)$.  Hence $\beta=W^*W$ commutes with 
the meridian $x$, so is a peripheral element and represents a slope
on the boundary torus of the knot exterior.

The exponents $\e(i)$ can be more explicitly described. There is an
odd integer $k$ coprime to $2n+1$ such that 
$$\e(i)=(-1)^{\lfloor \frac{ik}{2n+1}\rfloor}$$
 for each $i$.  In particular, since 
$$\frac{ik}{2n+1}-1<\lfloor \frac{ik}{2n+1}\rfloor< \frac{ik}{2n+1}$$
for each $i$, we have 
$$\lfloor \frac{ik}{2n+1} \rfloor + \lfloor \frac{(2n+1-i)k}{2n+1}\rfloor = k-1\equiv 0~\mathrm{mod}~2$$ for
each $i$, so that $\e(2n+1-i)=\e(i)$.  From this, it follows that
$$W(x^{-1},y^{-1})=W(y,x)^{-1}=W^{*-1}~~\mathrm{and}~~W^*(x^{-1},y^{-1})=W^{-1}.$$

\medskip
The following construction is essentially due to Burde (see
\cite[p.116]{B}).
Under the action of $SU(2)$ by rotations on $S^2$, we may choose 
fixed points of $\rho(x),\rho(W),\rho(y),\rho(W^*)$ as the vertices
$A,B,C,D$ respectively of a spherical rhombus, such that 
$\rho(W)(A)=C$ and $\rho(W^*)(C)=A$.  (There are degenerate
cases: possibly $A=C$ if $\rho(G)$ is abelian; possibly $B=D$
if $\rho(W)=\rho(W^*)$ with $\rho(W)^2=-I$.)  It follows that the
angle of rotation of $W^*W$ is $2\theta$ modulo $2\pi$, where $\theta$
is the angle $\widehat{DAB}$ of the rhombus.

\begin{center}
\epsfxsize=8cm \epsfysize=6cm \epsfbox{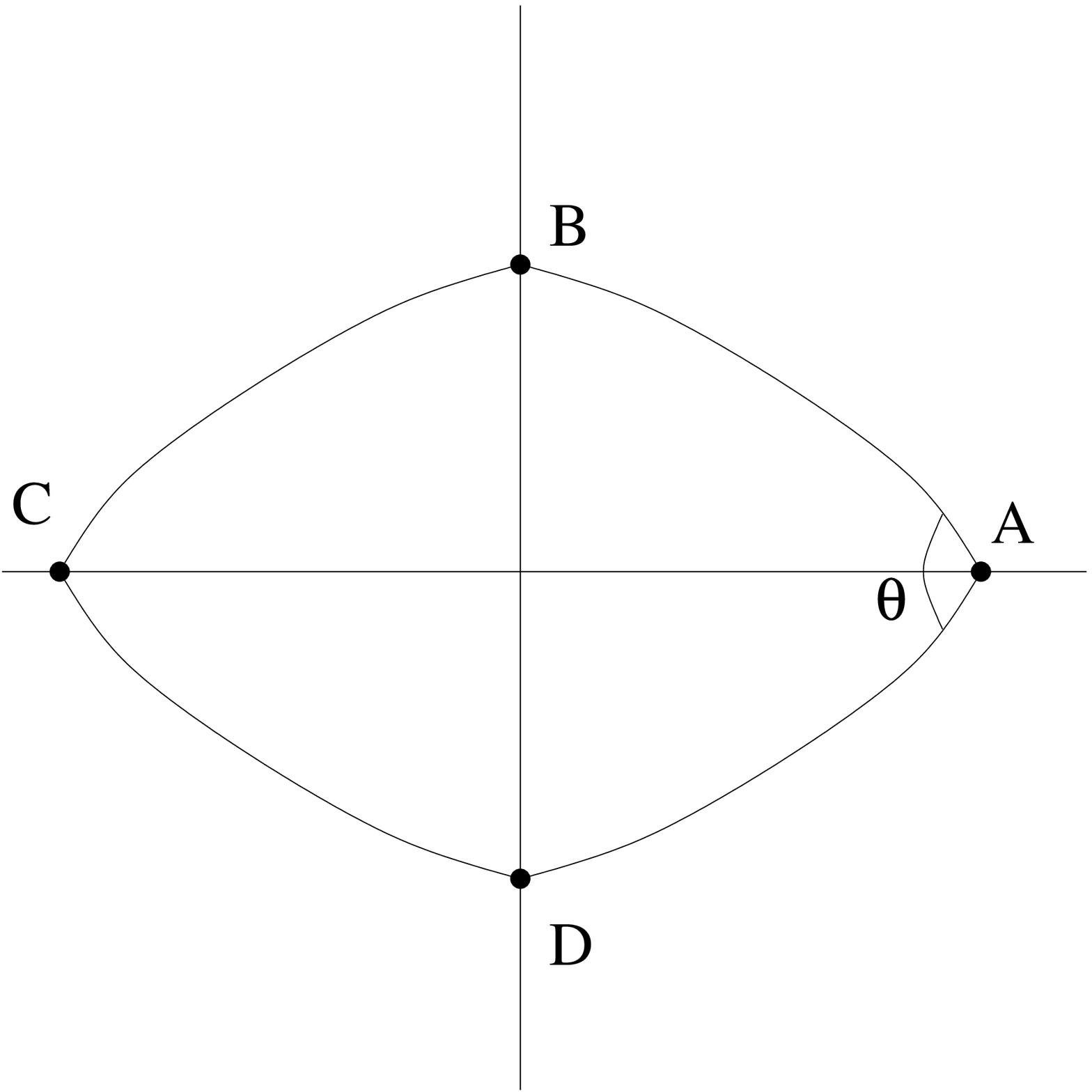}
\end{center}

Conjugacy in $SU(2)$ allows us freedom to place this rhombus where we
wish.  Let us choose to place it with $A=(1,0,0)$, and 
$C=(\cos\psi,\sin\psi,0)$ with $0<\psi<\pi$.  

If we have a path $\rho_t$ ($0\le t\le 1$) of representations, 
then this gives rise to a path $A_tB_tC_tD_t$ of rhombi, and a path
$\theta_t\in\R/(2\pi\Z)$ of corresponding angles.  Parameters $t$
with $\theta_t\in 2\pi\Z$ correspond to degenerate rhombi with $B_t=D_t$,
and hence to representations $\rho_t$ with $\rho_t(W)=\rho_t(W^*)$.

\medskip
Among all $SU(2)$ representations of $G$, a special r\^{o}le is played
by those whose image in $SO(3)$ is dihedral, in other words where
$\rho(x)^2=\rho(y)^2=-I$.  In this case, the points $B,D$
of our rhombus coincide with the north and south poles
$N,S=(0,0,\pm 1)$.   Burde \cite[pp. 116-117]{B}
explains that, if $\G$ is the group of a two-bridge knot which is not
a torus knot, then there is a path $\rho_t$ of irreducible representations
joining two dihedral representations $\rho_0,\rho_1$, such that $B,D$
switch poles on travelling from $t=0$ to $t=1$.  In other words, the
change of angle $\theta_1-\theta_0$ on traversing this path is an odd
multiple of $2\pi$ (in particular nonzero).  Replacing the path
$\rho_t$ by a smooth approximation if necessary, we may assume that
$\theta_t$ is differentiable as a function of $t$, and express this
as
$$\int_0^1 \frac{\partial\theta_t}{\partial t} dt \ne 0.$$

\medskip
Now consider another path of representations $\bar{\rho}_t$, defined
by $\bar{\rho}_t(x)=-\rho_t(x^{-1})$, $\bar{\rho}_t(y)=-\rho_t(y^{-1})$.
The equation $W(x^{-1},y^{-1})=W^{*-1}$ enables us to verify that
$\bar{\rho}_t$ is indeed a representation for each $t$.  Moreover,
since $\rho_t(x)^2=\rho_t(y)^2=-I$ for $t=0,1$, it follows that
$\bar{\rho}_t=\rho_t$ for $t=0,1$.   Finally, since $\bar{\rho}_t(W^*W)=\rho_t(W^*W)^{-1}$, the change in $\theta$ along
the path $\bar{\rho}$ is the negative of the change along the path 
$\rho_t$:
$$\frac{\partial\bar{\theta_t}}{\partial t}=-\frac{\partial\theta_t}{\partial t}.$$  
If $\gamma$ is the closed curve formed by concatenating
the paths $\rho_t$ and $\bar{\rho}_{1-t}$, the change in $\theta$ around
$\gamma$ is precisely twice that along $\rho_t$, namely an odd multiple of $4\pi$:
$$\int_\gamma \frac{\partial\theta_t}{\partial t} dt =
\int_0^1 \frac{\partial\theta_t}{\partial t} dt + \int_1^0 \frac{\partial\bar{\theta_t}}{\partial t} dt = 2 \int_0^1 \frac{\partial\theta_t}{\partial t} dt \ne 0.$$  
In particular $\nu([\gamma],W^*W)\ne 0$.
\end{proof}

\begin{cor}
Let $X$ be the exterior of a two-bridge knot in $S^3$, and let 
$X(\alpha)$ be the manifold formed from $X$ by Dehn filling along
a non-meridian slope $\alpha$ in $\partial X$.  Then $\pi_1(X(\alpha))$ 
admits an irreducible representation to $SU(2)$.
\end{cor}

\begin{proof}
By Theorem \ref{twobridgethm}, there is a closed curve $\gamma$ of irreducible
representations $pi_1(X)\to SU(2)$ such that the pairing $\nu$ on
$\pi_1(\gamma)\times\pi_1(\partial X)$ is not identically zero.

Then $\nu([\gamma],-):\pi_1(\partial X)\to\Z$ has kernel $\mu\Z$.  Since
$\alpha\notin\mu\Z$, $\nu([\gamma],\alpha)\ne 0$.  In other words,
the closed curve $t\mapsto\gamma_t(\alpha)\in S^1$ has non-zero winding 
number, and so is surjective.  There exists a point $\rho\in\gamma$ such that
$\rho(\alpha)=1$ in $SU(2)$.  Since $\pi_1(X(\alpha))$ is the quotient
of $\pi_1(X)$ by the normal closure of $\alpha$, $\rho$ induces
a representation $\pi_1(X(\alpha))\to SU(2)$ with nonabelian 
image.
\end{proof}

\section{Torus knots}\label{tk}

In this section we demonstrate that the pairing $\nu$ is not identically zero
on suitable curves in the $SU(2)$-representation
variety of a torus knot. We then apply this to 
the fundamental group of any manifold obtained
by nontrivial Dehn surgery on a torus knot, and study
its representations to $SU(2)$.

The $(p,q)$-torus knot has fundamental group $\Gamma=\<x,y|x^p=y^q\>$.
In particular, it has nontrivial centre, generated by $\zeta=x^p=y^q$.
If $\{\mu,\lambda\}$ is any meridian-longitude pair, then $\zeta$ belongs
to the peripheral subgroup generated by $\{\mu,\lambda\}$, since it commutes
with $\mu$.

The character variety $\mathcal{X}$
of $Hom(\Gamma,SU(2))$ splits into a number of
arcs as follows.  As for all knots, the subvariety
$\mathcal{X}_{red}$ corresponding to reducible representations
is isomorphic to the closed interval $[-2,2]$, parametrised by the
trace of $\rho(\mu)$.

If $\rho:\Gamma\to SU(2)$ is an irreducible representation, then
$\rho(x),\rho(y)$ are non-commuting matrices with $\rho(x)^p=\rho(y)^q$.
This can arise only if $\rho(x)^p=\rho(y)^q=\pm I$, where $I$
is the identity matrix.  Hence $\rho(x)$ has trace $2\cos(a\pi/p)$ and
$\rho(y)$ has trace $2\cos(b\pi/q)$ for some integers $a,b$ of the
same parity.  There are $(p-1)(q-1)/2$ open arcs $A_{(a,b)}$
in the irreducible character variety,
one corresponding to each pair $(a,b)$ of integers with $1\le a\le p-1$,
$1\le b\le q-1$, $a\equiv b$ modulo $2$.  Each open arc $A_{(a,b)}$
is the interior of a closed arc $\overline{A}_{(a,b)}$
in the whole character variety, whose endpoints
are reducible characters.

\begin{lemma}\label{endpoints}
The endpoints of $\overline{A}_{(a,b)}$ are the points
$$2\cos(c\pi/pq),2\cos(d\pi/pq)\in [-2,2]\cong\mathcal{X}_{red},$$
where where $c,d\in\{1,\dots,pq-1\}$
are the unique solutions to the congruences
$$c,d\equiv\pm a~~\mathrm{modulo}~2p;~~~c,d\equiv\pm b~~\mathrm{modulo}~2q.$$
\end{lemma}

\begin{proof}
On $A_{(a,b)}$, the trace of $\rho(x)$ is constant at $2\cos(a\pi/pq)$,
so the same will hold at each endpoint of $A_{(a,b)}$, which corresponds
to a reducible representation.  But $x\equiv\mu^{\pm q}$ modulo the
commutator subgroup, so for any reducible representation $\rho$ we
have $\rho(x)=\rho(\mu)^{\pm q}$.  If $z$ is a complex $q$-th root of $\cos(a\pi/pq)\pm i\sin(a\pi/pq)$, then $z=\cos(c\pi/pq)+i\sin(c\pi/pq)$
where $c\equiv\pm a~\mathrm{mod}~2p$.   Hence, for a reducible representation
$\rho$ at an endpoint of $A_{(a,b)}$, the trace of $\rho(\mu)$
must be $2\cos(c\pi/pq)$ with $c\equiv\pm a~\mathrm{mod}~2p$.

A similar analysis using $\rho(y)=\rho(\mu)^p$ gives the congruence
$c\equiv\pm b~\mathrm{mod}~2q$.

Finally, note that, since $a\equiv b~\mathrm{mod}~2$ and since $p,q$ are
coprime, each of the four pairs of simultaneous congruences
$$c\equiv\pm a~\mathrm{mod}~2p;\qquad c\equiv\pm b~\mathrm{mod}~2q$$
has a unique solution modulo $2pq$.  Moreover, if $c$ is the 
solution of one of these pairs of congruences, then $2pq-c$ is the 
solution of another, so precisely two of the four solutions lie in the
indicated range $\{1,\dots,pq-1\}$.
\end{proof}

\begin{prop}
Let $\gamma$ be the closed curve in $\mathcal{X}$
formed by the arc $\overline{A}_{(a,b)}$ together with
the subinterval $[2\cos(c\pi/pq),2\cos(d\pi/pq)]$
of $[-2,2]\cong\mathcal{X}_{red}$.  
Then $\nu([\gamma],\zeta)\ne 0$.
\end{prop}

\begin{proof}
The knot is embedded in an unknotted torus $T\subset S^3$.  Each
component of $S^3\sm T$ is an open solid torus.  Moreover, $x,y$
are represented by the cores of these solid tori, and $\zeta=x^p=y^q$
represents a curve on $T$ parallel to the knot.  In particular,
$\zeta\in A$, ie $\zeta$ is a peripheral curve.  Now $\rho(\zeta)=\pm I$ for
any irreducible representation $\rho$, and so $\rho(\zeta)$
is constant for $\rho\in A_{(a,b)}$.

Let $z=\exp(i\pi/pq)$, a primitive $(2pq)$-th root of unity.  Then the endpoints of $A_{(a,b)}$ correspond to the reducible representations
$\mu\mapsto z^c$ and $\mu\mapsto z^d$, where $c,d$ are given by Lemma
\ref{endpoints}.

Now, as $\rho$ moves continuously through reducible representations
from $\mu\mapsto z^c$ to $\mu\mapsto z^d$, the argument of 
$\rho(\mu)$ changes by $(d-c)\pi/pq$, so the argument of 
$\rho(\zeta)=\rho(\mu)^{pq}$ changes by $(d-c)\pi$, whence
$\nu([\gamma],\zeta)=(d-c)/2\ne 0$.
\end{proof}

\begin{cor}
Let $X$ be the exterior of a torus knot in $S^3$, and $X(\alpha)$
the manifold obtained from $X$ by Dehn filling along a non-meridian
slope $\alpha$.  Then $\pi_1(X(\alpha))$ admits a nontrivial
representation to $SU(2)$.
\end{cor}

\begin{proof}
If $\gamma$ is the curve in the Theorem, then $\nu([\gamma],\zeta)\ne 0$,
and so the kernel of the homomorphism $A\to\Z$, $\beta\mapsto\nu([\gamma],\beta)$, is precisely $\mu\Z$.  But by hypothesis
$\alpha\notin\mu\Z$, so $\nu([\gamma],\alpha)\ne 0$.  Thus the closed
curve $t\mapsto\gamma_t(\alpha)$ has nonzero winding number on $S^1$, so
is surjective.  There is a representation $\rho\in\gamma$ such that $\rho(\alpha)=1$ in $SU(2)$.  This choice of $\rho$ induces a nontrivial
representation $\pi_1(X(\alpha))\to SU(2)$.
\end{proof}

Of course, the above corollary is neither new nor surprising.  For example,
almost all the groups $\pi_1(X(\alpha))$ have nontrivial abelianisation, so
admit representations to $SU(2)$ that are reducible but nontrivial.
Of more interest is the question of which $\pi_1(X(\alpha))$ admit 
irreducible representations to $SU(2)$.  This question can also be readily
answered using the known classification of $3$-manifolds obtained by
Dehn surgery on torus knots \cite{M}.  Here we present an alternative approach
using an adaptation of our winding-number technique.

\begin{thm}
Let $X$ be the exterior of the $(p,q)$ torus knot, where $1<p<q$,
and $X(\alpha)$
the manifold obtained from $X$ by Dehn filling along a non-meridian
slope $\alpha\in\Q\cup\{\infty\}$.  Then
\begin{enumerate}
\item[\rm{(1)}] if $\alpha=pq$ and $p>2$, then $\pi_1(X(\alpha))$ admits
an irreducible representation to $SU(2)$;
\item[\rm{(2)}] if $\alpha=pq$ and $p=2$, then $\pi_1(X(\alpha))$ admits
no irreducible representation to $SU(2)$, but admits
a representation to $SO(3)$ with nonabelian image;
\item[\rm{(3)}] if $\alpha=pq\pm\frac1n$ for some positive integer
$n$, then every representation from $\pi_1(X(\alpha))$ to $SO(3)$
has abelian image;
\item[\rm{(4)}] for any other value of $\alpha$, $\pi_1(X(\alpha))$ admits
an irreducible representation to $SU(2)$.
\end{enumerate}
\end{thm}

\noindent{\bf Remark} The statement of this theorem fits the classification
of \cite{M}, where it is proved that $X(\alpha)$ is a lens
space in Case (3); a connected sum of two lens spaces in Cases
(1) and (2);
and a Seifert fibre space in Case (4).

\medskip
\begin{proof}

(1) Since $2<p<q$, one of the components of $\mathcal{X}_{irr}$ is the arc
$A_{(2,2)}$.  But any point on $A_{(2,2)}$ corresponds to a representation
$\rho$ with $\rho(x^p)=\rho(y^q)=I$.

\medskip
(2) In this case $\pi_1(X(\alpha))\cong\Z_2\ast\Z_q$.  Since the only 
element of order $2$ in $SU(2)$ is the central element $-I$, the image of
any representation $\Z_2\ast\Z_q\to SU(2)$ is abelian.  However, corresponding
to any point on $A_{(1,1)}$ is a representation $\rho$ with
$\rho(x^2)=\rho(y^q)=-I$, so composing this with the quotient map
$SU(2)\to SO(3)$ gives a representation of $\pi_1(X(\alpha))$ to $SO(3)$
with nonabelian image.

\medskip
(3) Let $\zeta$ be the curve $x^p=x^q$ of slope $pq$.  Then $\zeta=\mu^{pq}\lambda$, so $\alpha=\mu^{npq\pm 1}\lambda^n=\mu^{\pm 1}\zeta^n$
in $\pi_1(\partial X)$.  Now any representation from $\pi_1(X(\alpha))$
to $SO(3)$ with nonabelian image arises from a representation of
$\pi_1(X)$ with nonabelian image, which therefore lifts to an irreducible
representation $\rho:\pi_1(X)\to SU(2)$, such that $\rho(\alpha)=\pm I$.
But $\rho$ corresponds to a point on one of the open arcs $A_{(a,b)}$,
so $\rho(\zeta)=(-I)^a$ and hence $\rho(\mu)=(\rho(\alpha)\rho(\zeta)^{-n})^{\pm 1}=\pm I$, contradicting the assumption that $\rho$ is irreducible.

\medskip
(4) As in the previous case, let $\zeta=\mu^{pq}\lambda$ denote the curve
with slope $pq$.  Then $\pi_1(\partial X)$ is generated by $\zeta$ and $\mu$,
so we can write $\alpha=\mu^g\zeta^h$.  If $|g|\le 1$ then we are in one
of the previous cases, so we have $|g|\ge 2$.

Suppose first that $pq$ is even.   Then the endpoints of $A_{1,1}$ are reducible
representations $\rho$ in which the trace of $\rho(\mu)$ is $\pm 2\cos(\pi/pq)$.
Choose $\theta\in [\pi/pq,(pq-1)\pi/pq)]$ such that $\theta$ is an odd
multiple of $\pi/|g|$.  Then by continuity of trace, we can choose $\rho\in
A_{(1,1)}$ such that the trace of $\rho(\mu)$ is $2\cos(\theta)$.  Provided
$h$ is odd, this gives $\rho(\mu)^g=-I=\rho(\zeta)^{-h}$, so $\rho(\alpha)=I$.
If $h$ is even then $|g|$ is odd, since $\alpha$ is a slope.
In particular $|g|>2$. In this case, we take $\theta$ to be an even multiple 
of $\pi/|g|$, and the argument goes through as before.

Now consider the case where $pq$ is odd. Precisely one of the two
positive integers $(q\pm p)/2$ is odd.  Call it $c$, and note that $c\in\{1,\dots,q-1\}$.  Let $a$ be the unique odd integer with
$1\le a\le p-1$ and $a\equiv\pm~c~\mathrm{mod}~p$.
Then the endpoints of $A_{a,c}$ are reducible representations $\rho$ where
the trace of $\rho(\mu)$ is $2\cos(c\pi/pq)$ and $2\cos((pq-q+c)\pi/pq)$
respectively.  Now the interval $[c\pi/pq,(pq-q+c)\pi/pq]$ contains at least
one odd multiple of $\pi/|g|$, and (if $|g|>2$) at least one even
multiple of $\pi/|g|$.   Arguing as before, we can choose $\rho\in A_{a,c}$
such that $\rho(\mu)^g=\rho(\zeta)^{-h}$, and so $\rho(\alpha)=I$, except 
possibly if $|g|=2$ and $h$ is even (which does not arise, since
$\alpha$ is a slope).

\end{proof}

\end{document}